\documentclass[12pt]{article}
\usepackage[latin1]{inputenc}
\usepackage{amssymb}
\usepackage{latexsym}
\usepackage{amscd}
\usepackage{longtable}
\usepackage{amsmath}
\usepackage{array}
\usepackage[all]{xy}
\usepackage{a4}
\usepackage{mathptmx}
\newtheorem{teo}{Theorem}[section]
\newtheorem{lema}{Lemma}[section]
\newtheorem{coro}{Corollary}[section]
\newtheorem{prop}{Proposition}[section]

\newtheorem{rem}{Remark}[section]

\newcommand{\Q}{\mathbb Q}

\newcommand{\HH}{\mathbb H}
\newcommand{\R}{\mathbb R}
\newcommand{\Z}{\mathbb Z}

\newcommand{\K}{\mathbb K}
\newcommand{\F}{\mathbb F}
\newcommand{\E}{\mathbb E}
\newcommand{\M}{\mathbb M}
\newcommand{\LL}{\mathbb L}

\newcommand{\C}{\mathbb C}
\newcommand{\GL}{\mathrm{GL}}

\newcommand{\End}{\operatorname{End}}
\newcommand{\Frob}{\operatorname{Frob}}
\newcommand{\ord}{\operatorname{ord}}

\newcommand{\rank}{\operatorname{rank}}

\newcommand{\Gal}{\operatorname{Gal}}

\newcommand{\Res}{\operatorname{Res}}

\newcommand{\Norm}{\operatorname{N}}

\newcommand{\cO}{\mathcal {O}}

\newcommand{\cS}{\mathcal {S}}

\newcommand{\cI}{\mathcal {I}}
\newcommand{\cP}{\mathcal {P}}

\newcommand{\cK}{\mathcal {K}}
\newfont{\gotip}{eufb10 at 12pt}
\newcommand{\gn}{\mbox{\gotip n}}
\newcommand{\gm}{\mbox{\gotip m}}

\newcommand{\ga}{\mbox{\gotip a}}

\newcommand{\gp}{\mbox{\gotip p}}

\newcommand{\li}[1]{\text{li}\left(#1\right)}

\title
{Cropping Euler factors of modular $L$-functions}
\author{Josep Gonz\'{a}lez, Jorge Jim\'enez-Urroz, Joan-Carles Lario
\footnote{The first author is partially
supported  by DGICYT Grant MTM2009-13060-C02-02, the second by DGICYT
Grant MTM2009-11068 and the third by DGICYT Grant MTM2009-13060-C02-01}}

\date{\today}

\begin{document}

\maketitle

\begin{abstract}
As it is well-known, much of the arithmetic information for a
Galois number field extension $\LL / \Q$ is encoded by its
Dedekind zeta function and the set of primes that split completely
in $\LL$. According to the Birch and Swinnerton-Dyer conjectures,
if $A / \Q$ is an abelian variety then its $L$-function must also
capture a substantial part of the properties of~$A$. The smallest
number field $\LL$ where $A$ has all its endomorphisms defined
must also have a role. This article deals with the relationship
between these two objects in the specific case of modular abelian
varieties $A_f / \Q$ associated to weight 2 newforms for the group
$\Gamma_1 (N)$. Specifically, our goal is to relate the order of
$L (A_f/\Q, s)$ at $s = 1$ with Euler products cropped by
primes that split completely in $\LL$. This is attained by giving separated formulae for the CM and non CM cases when a power of $A_f$ is isogenous over $\Q$ to the Weil restriction of the building block of $A_f$.
\end{abstract}

\section{Introduction}
Let $f$ be a normalized modular newform in $S_2(\Gamma_1(N))$
with Fourier expansion given by $\sum_{n>0} a_n q^n$ and let
$\varepsilon$ be its Nebentypus. We shall be concerned with its
$L$-function
$$
L(f,s)=\sum_{n>0}\frac{a_n}{n^s}=\prod_{p}\frac{1}{1-a_p p^{-s}+ \varepsilon (p) p^{1-2s}}\,,
$$
and, specially, with its Euler product. The function $L(f,s)$
converges absolutely  for $\Re(s)>3/2$, and has analytic
continuation to the whole complex plane. It is known that
$\E=\Q(\{ a_n\})$ is a number field and the Galois action on the
Fourier coefficients provides a set of normalized newforms
$f_1,\cdots,f_n$ of $S_2(\Gamma_1(N))$ with $n=[\E:\Q]$. The
product $\prod_{i=1}^n L(f_i,s)$ is the $L$-function $L(A_f/\Q,s)$
of the $n$-dimensional abelian variety $A_f/\Q$ attached by
Shimura to $f$. The value $\ord_{s=1} L(A_f/\Q,s)$ is a matter of
importance since it must coincide with the rank of the
Mordell-Weil group $A_f(\Q)$ according to the Birch and
Swinnerton-Dyer conjectures.

To motivate the issue that we want to address in this article, let
us consider for a moment the case of Euler products arising from
Dedekind zeta-functions attached to number fields. For every
number field $\LL$,  it is well-known that its Dedekind
zeta-function $\zeta_{\LL}(s)$, defined as an Euler product on the
right half-plane $\Re(s)>1$,  does have meromorphic continuation
to the whole complex plane and it has a unique simple pole at
$s = 1$. For the particular case that $\LL/\Q$ is a Galois
extension, we consider the subset $\cS_1$ of rational primes that
split completely in $\LL$ and introduce the partial Dedekind
zeta-function
$$ \displaystyle{\zeta_\LL(\cS_1, s):=\prod_{p\in\cS_1}\frac{1}{1-p^{-s}}}\,,
$$
defined on  $\Re(s)>1$. Since $\zeta_{\LL}(s)/\zeta_\LL(\cS_1,
s)^{[\LL:\Q]}$ is holomorphic on $\Re(s)>1/2$ and does not
vanishes at $s=1$, the function  $\zeta_\LL(\cS_1, s)^{[\LL:\Q]}$
admits meromorphic continuation on $\Re(s)>1/2$ and satisfies
$$ \ord_{s=1} \zeta_\LL(\cS_1, s)^{[\LL:\Q]}=\ord_{s=1
}\zeta_{\LL}(s)=-1\,.
$$
We point out  that this equality shows that $\zeta_\LL(\cS_1, s)$
does not admit meromorphic continuation on $\Re(s)>a$ for any
$a<1$, except for the trivial case $\LL=\Q$.

The starting point of this article is to study the generalization
of this phenomenon  to modular $L$-functions. In other words, we
want to find out if the Euler product of $L(f, s)$ can be cropped
in the  sense that it exists a subset of distinguished primes with
regard to their contribution to the order of $L (f, s)$ at $s =1$.

It turns out that there is a natural place to look at for finding
this set of primes. Indeed, the splitting field $\LL$ of $A_f$
(that is, the smallest number field where $A_f$ has all its
endomorphisms defined) is an important ingredient of the
arithmetic of $A_f$. In particular, the abelian variety $A_f$ is
isogenous over $\LL$ to the power of a simple  abelian variety
$B_f$. In \cite{gola01}, \cite{gola07} and \cite{go09}, the field
$\LL$ is explicitly determined. Then, we propose to consider the
partial Euler product
$$
L(f,\cS_1, s) := \prod_{p\in \cS_1} \frac{1}{1-a_p p^{-s}+ \varepsilon (p)
p^{1-2s}}\,,
$$
where $p$ runs over the set $\cS_1$ of primes that split
completely in $\LL$.

The plan of this paper is as follows. Section 2 is  devoted to
introduce notation and summarizes some well-known facts concerning
modular $L$-functions. Since two different situations emerge
depending on whether $f$ has complex multiplication (CM) or not,
each of them is treated separately in  Sections 3 and 4. In both
sections, we study the relationship between
$\ord_{s=1}L(A_f/\LL,s)$ and the order of
$\prod_{i=1}^nL(f_i,\cS_1,s)^{[\LL:\Q]}$ at $s=1$. As a good
point, the function $L(f,\cS_1, s)$ does not depend on certain
Galois conjugates of the newform~$f$ but, as we shall see,
unfortunately the primes $p$ with residue degree $d(p)=2$ in $\LL$
(if any) will cause some problems and we shall need a substitute
of $L(f,\cS_1, s)$ as first approach. In both cases, we introduce a partial
Euler product $L(s)$ of $L(A_f/\Q,s)$ associated with primes $p$
with residue degree $d(p)\leq 2$ in~$\LL$. We prove in Theorem~3.1
and Theorem~4.1 that $\ord_{s=1}L(B_f/\LL,s)$ is determined by
$\ord_{s=1}L(s)^{[\LL:\Q]/2}$. For the particular case that a
power of $A_f$  is isogenous over $\Q$ to the  Weil restriction
$\Res_{\LL/\Q}(B_f)$, then $L(B_f/\LL,s)$  agrees with the
corresponding power of  $L(A_f/\Q,s)$ and we obtain results for
this $L$-function.

Section~5 contains the main results of the article. This last
section is devoted to study whether we can avoid the primes $p$
with $d(p)=2$ in order to use the more natural $L(f,\cS_1, s)$
instead of $L(s)$. As we will show, this fact is related to the
distribution of the values $b_p=(a_p^2-2\varepsilon(p))/(2\,p)$
for primes~$p$ with~$d(p)=2$. For the CM case, we generalize
results of T. Mitsui in~\cite{mitsui} on distribution of primes in
sectors and this allows us to solve completely this problem in
Theorem~\ref{mainCM}. For the non-CM case,
we extend the recent result obtained by T. Barnet-Lamb, D. Geraghty, M.
Harris and R. Taylor in \cite{BGHT} about Sato-Tate distributions, when
restricting to primes in arithmetic progressions.
Then, we present the main result
for this case in Theorem~\ref{mainnoCM}.

\section{Modular $L$-functions}
Let $f=\sum_{n>0}a_n q^n$ be a normalized newform of level $N$
with Nebentypus $\varepsilon$ and let~$\E$ be the number field
$\Q(\{a_n\})$. From now on, at our convenience an $L$-functions with
an asterisk will stand for the corresponding $L$-functions but
removing the Euler factors attached to the primes dividing $N$.

Let $\lambda$ be a prime ideal of $\E$ over a rational prime $\ell$.
There is a continuous $\lambda$-adic representation
$$\rho_\lambda:\Gal(\overline{\Q}/\Q) \longrightarrow \GL_2(\cO_\lambda)\,,$$
where $\cO_\lambda$ denotes the completion of the ring of integers
of $\E$ at $\lambda$ such that $L(f,s)$ is the $L$-function attached
to this $\lambda$-adic representation. Also, for every number
field~$\LL$, we shall denote by $L(f/\LL,s)$ the $L$-series
attached to $\rho_\lambda$ restricted to $\Gal
(\overline{\Q}/\LL)$. If $\LL/\Q$ is a Galois extension and $d(p)$
denotes the residue degree of $p$ in~$\LL$, then one has
$$
L^*(f/\LL,s)=\prod_{p\in\cP}\frac{1}{(1- b_p p^{- d(p)\,s}+ \varepsilon (p)^{d(p)}
p^{d(p)(1-2s)})^{[\LL:\Q]/d(p)}}\,,
$$
where hereafter $\cP$ will denote the set of rational primes not
dividing the level $N$,
$$b_p={\alpha_p}^{d(p)}+(\overline{\alpha}_p \varepsilon(p))^{d(p)}\,,$$
and $\alpha_p$ is any root of the polynomial $x^2-a_p\,x+
p\,\varepsilon(p)$. Notice that if $p\neq \ell$, then
$$
b_p = \operatorname{Trace} \left( \rho_\lambda(\Frob_p )^{d(p)}
\right)\,.
$$
By the results on  base change of automorphic representations
established by Langlands~\cite{Langlands}, we know that $ L^*(
f/\LL,s)$ has analytic continuation to the whole complex plane
when the group $\Gal (\LL/\Q)$ is solvable and, consequently,  so
does
 $$L^*(A_ f/\LL,s)=\prod_{\sigma: \E\hookrightarrow
\overline{\Q}} L^*({}^\sigma f/\LL,s)\,.$$
Moreover, it is clear that $L^*(f/\LL,s)$ and $L(f/\LL,s)$
have the same order at $s=1$. From now on, the number field $\LL$
will be the splitting field of $A_f$; i.e., $\LL$ is the smallest
number field where $A_f$ has all its endomorphisms defined. Hence,
$A_f$ is isogenous over~$\LL$ to the power of an absolutely simple
abelian variety $B_f/\LL$, the so-called building block of~$A_f$, and one has
$$
L^*(A_f/\LL,s)=L^*(B_f/\LL ,s)^{\dim A_f/\dim B_f}\,.
$$
Moreover, by Milne \cite{milne72}, one has
$$
L^*(B_f/\LL ,s)=L^*(\Res_{\LL/\Q}(B_f),s)\,.
$$
As we shall show, the function $L^*(B_f/\LL ,s)$ is the product of
functions  $ L^*({}^\sigma f/\LL,s)$ when $\sigma$ runs over a
certain subset of the embeddings of $\E$ into $\overline\Q$ and,
thus, it has also analytic continuation to the whole complex
plane since $\Gal(\LL/\Q)$ is solvable.

To end this section, let us fix the following terminology. For every
subset $ \cS \subseteq\cP$, we consider the partial Euler products
$$\begin{array}{ll}
L(f^{\phantom{-}},\cS,s):=&\displaystyle{\prod_{p\in\cS} \frac{ 1}{1-a_p
p^{-s}+\varepsilon(p) p^{-2s+1}}}\,,\\[6 pt]
 L(f^{-},\cS,s):=& \displaystyle{\prod_{p\in\cS} \frac{ 1}{1+a_p
p^{-s}+\varepsilon(p) p^{-2s+1}}\,,}
\end{array}
$$
where we take $L(f,\emptyset,s)=L(f^{-},\emptyset,s)=1$. One
has $L^*(A_f/\Q,s)=\prod_{\sigma} L({}^\sigma f,\cP,s)$, where
$\sigma$ runs over the set of embeddings of $\E$ into
$\overline{\Q}$.

\section{First approach to the non-CM case}
In this section we assume that $f$ is without CM. We recall that a
Dirichlet character $\chi$ is called an inner-twist of $f$
if there is an embedding
$\sigma: \E \hookrightarrow \overline{\Q}$ satisfying
${}^{\sigma}a_p= \chi (p) a_p$ for all primes  $p\in \cP$. If for
an embedding $\sigma$ there is an inner-twist, it is unique and is
denoted by  $\chi_\sigma$.

By Proposition 2.1 in \cite{gola01}, the splitting field $\LL$ of
$A_f$ is the number field $\overline{\Q}^{\cap_{\sigma}\ker
\chi_{\sigma}}$, where $\sigma$ runs over the set of embeddings of
$\E$ into $\overline{\Q}$ for which there is  an inner-twist
$\chi_\sigma$. The extension $\LL/\Q$  is the compositum of the
cyclic extension $\overline{\Q}^{\ker \varepsilon}$ and a
polyquadratic extension of $\Q$.  Notice that  $\LL$ is contained
in the $N$-th cyclotomic field and, thus, all primes in $\cP$ are
unramified in $\LL$. Moreover, the center of the algebra
$\End_{\LL} (A_f)\otimes\Q $ is the totaly real subfield $\F=\Q(\{
a_p^2/\varepsilon(p): p \in \cP \})$ of $\E$ and $\dim B_f= t
\cdot [\F :\Q]$, where $t$ is either $1$ or $2$ depending on
whether the algebra $\End_{\LL}(B_f)\otimes \Q$ is isomorphic to
either $\F$ or a quaternion algebra with center $\F$. In
particular, $A_f$ is isogenous over $\LL$ to $B_f^{[\E:\F]/t}$. We
shall need the following result.

\begin{lema}\label{tes1}
The abelian variety $A_f$ is isogenous over $\Q$
to $\Res_{\LL/\Q}(B_f)$ if and only if  $t=1$ and $[\LL:\Q]=[\E:\F]$.
\end{lema}
\noindent{\bf Proof.} Since $A_f$ is simple over $\Q$, $A_f$ is
isogenous over $\Q$ to $\Res_{\LL/\Q}(B_f)$ if and only if
$[\E:\Q]=[\LL:\Q] \dim B_f$, i.e. $[\E:\F]=[\LL:\Q] \,t$. We know
that $\E/\F$ is an abelian extension and $\Gal(\E/\F)$ is the set
of embeddings of $\E$ into $\overline{\Q}$ for which there is  an
inner-twist of $f$. Therefore, $[\E:\F]\leq [\LL:\Q]$, and it
follows the lemma. \hfill $\Box$

\vskip 0.1 cm

As in Section~2, for a prime $p\in\cP$, let  $d(p)$ be the residue degree of $p$ in $\LL$
and $b_p$ be the trace of $\rho_\lambda(\Frob_p )^{d(p)}$.
Since $\overline{\Q}^{\ker \varepsilon}\subseteq \LL$, one has
$\varepsilon(p)^{d(p)}=1$. By Proposition~5.2 and Lemma~6.1
of~\cite{bago}, we know that for almost all primes such that $a_p\neq
0$, then $a_p^{d(p)}\in\F$ and $d(p)$ is the smallest positive
integer satisfying this condition. Moreover, $b_p\in \F$ for all
primes $p\in\cP$. We shall consider the following partition of~$\cP$:
\begin{equation}\label{s1}
\begin{array}{ll}
\cS_1=& \{ p\in\cP: d(p)=1\}\,,\\ [4 pt]
 \cS_2=&\{p\in\cP:
d(p)=2\}\,,\\[4 pt]
\cS_3=&\{ p\in\cP: d(p)\geq 3\}\,.
\end{array}
\end{equation}

Notice that for every  prime $p\in\cS_2$ there exists
$\tau\in\Gal(\E/\F)$ such that ${}^{\tau}a_p=-a_p$ and, thus,
$L(f,\cS_2,s)\,L(f^{-},\cS_2,s)$ is a partial Euler product of
$L(A_f/\Q,s)$.

\begin{teo}\label{teoremanoCM} Keep the above notations. Let $\cI$ be a minimal
set of embeddings of $\E$ into $\overline{\Q}$ such that their
restrictions on $\F$ provide all embeddings of $\F$ into
$\overline{\Q}$. Then,
\begin{itemize}
\item[(i)] the function
$$\left(L(f,\cS_1,s)^2\,L(f,\cS_2,s)L(f^{-},\cS_2,s)\right)^{[\LL:\Q]/2}$$
has analytic continuation to the right half-plane $\Re (s)>5/6$ and its order at
$s=1$ is equal to $\ord_{s=1} L(f/\LL,s)$.

\item[(ii)]
Let
$$
L_1(s):=\prod_{\sigma\in\cI} L({}^{\sigma}f,\cS_1,s)\,,\quad L_2(s):=
\prod_{\sigma\in\cI}L({}^{\sigma}f,\cS_2,s)\, L({}^{\sigma}f^{-},\cS_2,s)
$$
and $L(s):=L_1(s)^2 L_2(s)$. The function
$L( s)^{ [\LL :\Q]/2}$, that for $\cS_2=\emptyset$ coincides with  $L_1(s)^{ [\LL :\Q]}$,
satisfies
$$
\begin{array}{lr}
\ord_{s=1} L(B_f/\LL,s)=&\displaystyle{ t\cdot \ord_{s=1} L(s)^{ [\LL:\Q]/2}}, \\[5pt]
\ord_{s=1} L(A_f/\LL,s)=&\displaystyle{ [\E:\F]\cdot \ord_{s=1} L(s)^{ [\LL:\Q]/2}}
\,.
\end{array}$$

\item[(iii)]In the particular case that $A_f$ is isogenous to
$\Res_{\LL/\Q} (B_f)$,  then
$$
\ord_{s=1} L(A_f/\Q,s)=\displaystyle{\ord_{s=1}
 L(s)}^{[\E :\F]/2} \,.
$$
\end{itemize}
\end{teo}
\noindent {\bf Proof.} We consider the factorization  $L^*
(f/\LL,s)=\prod_{i=1}^3 G_i(f,s)$, where
$$
\begin{array}{ll}
G_1(f,s):=&\displaystyle{L(f,\cS_1,s) ^{[\LL:\Q]}},\\[4 pt]
G_2(f,s):=&\displaystyle{\left(\prod_{p\in\cS_2}\frac{1}{ 1-b_p p^{- 2\,s}+
 p^{2 (-2 s+1)}}\right)^{[\LL:\Q]/2}},\\[4 pt]
G_3(f,s):=&\displaystyle{\left(\prod_{p\in
\cS_3}\frac{1}{ 1-b_p p^{- d(p)\,s}+ p^{d(p) (-2s+1)}}\right)^{[\LL:\Q]/d(p)}}.
\end{array}
$$
Observe that if $p\in\cS_2$, then we have
$$
1- b_p p^{-2s}+p^{2(1-2 s)}= (1- a_p p^{-s}+\varepsilon(p)p^{1-2 s})(1+  a_p p^{-s}+\varepsilon(p)p^{1-2 s})\,.
$$
Therefore, it follows
\begin{equation}\label{eq1}
 L^*(f/\LL,s) = \left(L(f,\cS_1,s)^2\,L(f,\cS_2,s)L(f^{-},\cS_2,s)\right)^{[\LL:\Q]/2} \cdot G_3(f,s)\,.
\end{equation}
Due to the fact that $L^*(f/\LL,s)$ has analytic continuation to
the whole complex plane and (\ref{eq1}), in order to prove part
(i) it is enough to prove that $G_3(f,s)$ is analytic on  $\Re
(s)>5/6$ and $G_3(f,1)$ is non-zero. However, this follows from
$\mid  b_p\mid\leq 2 p^{d(p)/2}$ and the inequality
$$\mid b(p)
p^{-s d(p)}-p^{d(p) (-2 s+1)}\mid\leq 2 p^{-d(p)(\Re(s)-1/2)}+
p^{-2d(p)(\Re(s)-1/2)}\leq 3 p^{-3(\Re(s)-1/2)}\,,$$
valid since $d(p)\geq 3$.
Now, observe that for every $\sigma\in\Gal(\E/\F)$, we have
$G_i(f,s)=G_i({}^\sigma f,s)$ for all $i\leq 3$. Hence, it follows
part (ii). Finally, assume that $A_f$ is isogenous over  $\Q$ to
$\Res_{\LL/\Q}(B_f)$. By Lemma \ref{tes1}, we have that  $t=1$ and
$[\LL:\Q]=[\E:\F]$. Noting that
$L^*(A_f/\Q,s)=L^*(B_f/\LL,s)$, the last assertion of the
statement follows.
 \hfill $\Box$

\begin{rem}\label{rem2}
The above theorem has full sense when $[\LL:\Q]>1$. In this case,
if $\cS_2=\emptyset$ then $L({}^\sigma f,\cS_1,s)^{[\LL:\Q]}$ has
analytic continuation on $\Re(s)>5/6$ for all $\sigma\in\cI$ and,
moreover,
\begin{equation}\label{S1}  \ord_{s=1} L({}^\sigma f/\LL,s)=
\ord_{s=1} L({}^\sigma f,\cS_1,s)^{[\LL:\Q]} \,. \end{equation}
For $\cS_2\neq \emptyset$, the product
 $$
G_\sigma(s):=L({}^\sigma
f,\cS_2,s)L({}^{\sigma}f^{-},\cS_2,s)=\prod_{p\in\cS_2}\frac{1}{
1-{}^{\sigma}b_p p^{- 2\,s}+
 p^{2 (-2 s+1)}}
$$
converges absolutely for $\Re (s)>1$ since for $p\in\cS_2$ we have
$\mid {}^{\sigma}b_p \mid \leq 2 p$. In particular, $L({}^\sigma
f,\cS_1,s)^{[\LL:\Q]}$ is analytic  on $\Re(s)>1$.
As we shall show in the last section, the product $G_\sigma$
converges at $s=1$ if and only if the series
$$
\sum_{p\in\cS_2} \frac{{}^\sigma b_p}{p^2}
$$
converges.  In this case, the equality (\ref{S1}) also applies,
understanding, here and in the sequel, that if for a meromorphic function $H(s)$ defined on $\Re (s)>1$ there exists an integer $n$ satisfying that
$$
\lim_{s\to 1,\, \Re(s)>1} \frac{H(s)}{(s-1)^n}
$$
is a non-zero complex number, then  we write $\ord_{s=1} H(s)=n$.
\end{rem}

\section{First approach to the CM case}\label{CM}
Let $\K$ be an imaginary quadratic field and let $\cO$ be its ring
of integers. Now, assume that  $f=\sum_{n>0} a_n q^n$ has CM  by
$\K$.  Thus, there exist an integral ideal $\gm$ of $\K$ and a
primitive Hecke character $\psi: I (\gm)\to \C^*$ of conductor
$\gm$ such that $f=\sum  \psi(\ga)q^{\Norm (\ga)}$, where the
summation is restricted to the integral ideals $\ga$ of $\K$
coprime to $\gm$. Here,   $\Norm(\ga)$ is the norm of the ideal
$\ga$ and $I(\gm)$ denotes the multiplicative group of fractional
ideals of $\K$ relatively prime to $\gm$. In this case, the
$L$-function attached to $f$ can be rewritten as
$$
L(f,s)=\sum_{(\ga,\gm)=1} \frac{\psi (\ga)} {\Norm (\ga)^s}=\prod_{(\gp,\gm)=1}\frac{1}{1-\psi(\gp)
\Norm (\gp)^{-s}}\,,
$$
and the level $N$ of $f$ is $\Norm (\gm)$ times the absolute value
of the discriminant of $\K$. Note that for a prime $p\in\cP$, the
roots of the polynomial $x^2-a_p\,x+p\,\varepsilon(p)$ are
$\psi(\gp)$  and $\psi(\overline{\gp})$ when $p$ splits in $\K$
and  $\pm\sqrt{\psi (\gp)}$ for $p$ inert, where $\gp$ is a prime
of $\K$ over~$p$.

Attached to $\psi$  there is a character $\eta:(\cO/\gm)^*\to\C^*$
defined by $\eta(a)=\psi (a\cO)/a$. The Nebentypus $\varepsilon$
of $f$ is the Dirichlet character mod $N$ such that
$\varepsilon(n)=\eta(n)\chi(n)$, where $\chi$ denotes the
quadratic Dirichlet character attached to $\K$. The existence of
$\psi$ implies that the natural projection  $\cO^* \to
(\cO/\gm)^*$ is a group monomorphism and, thus, two different
generators of a principal ideal in $I(\gm)$ are not equivalent mod
$\gm$ and the unique unity of $\K$ in $\ker \eta$ is $1$.

\vskip 0.1 cm

We fix the following notation. For a subset  $S$ of $(\cO/\gm)^*$,
we denote by $P_S(\gm)$ the subset of $I(\gm)$ consisting on
principals ideals which have a generator $\alpha$ such that
$\alpha\pmod {\gm}$  lies in $S$. To simplify notation, we write
$P_{\delta}(\gm)$ and $P(\gm)$ when $S=\{\delta\}$ and
$S=(\cO/\gm)^*$, respectively. Of course, if $G$ is a subgroup of
$(\cO/\gm)^*$, then $P_G(\gm)$ is also a subgroup of $I(\gm)$
containing the subgroup $P_1(\gm)$.

Let $\K_{\gm}$ be the ray class field mod $\gm$. By Class field
theory we know that the group $I(\gm)/P_1(\gm)$ is isomorphic to
$\Gal (\K_{\gm}/\K)$ via the Artin map and there exists an
intermediate number field $\LL_{\eta}$ between the Hilbert class
field of  $\K$ and $\K_{\gm}$ such that $\Gal
(\LL_{\eta}/\K)\simeq I(\gm)/P_{\ker \eta}(\gm)$.
In~\cite{gola07}, it is proved that there is a quotient abelian
variety $A$ of $A_f$ defined over $\K$, simple over $\K$ and such
that $ A_f/\K$ is isogenous over~$\K$ to either $A\times \overline
A$ or $A$ according to whether  $\K$ is contained in $\E$ or not,
where $\overline{\phantom{c}}$ stands for the complex conjugation.
Moreover, the splitting field of $A$ is the number field
$\LL_\eta$, which is a cyclic extension of the Hilbert class field
of $\K$ (cf. Theorem~1.2 of \cite{gola07}). In Remark 2.1 of
\cite{go09}, it is showed that the splitting field of $A_f$ is the
compositum
$$
\LL=\LL_{\eta}\, \overline{\Q}^{\ker
\varepsilon}\,.
$$

Notice that again $\LL/\Q$ is solvable. Since $\Res_ {\K/\Q} (A)$
is isogenous over $\Q$  to either $A_f$ or $A_f^2$ depending on
whether $\K\subseteq \E$ or not, we have that
$$
L(A/\K,s)= L(A_f/\Q,s)^{[\E\,\K: \E]}\,.
$$
Now, the building block of $A_f$ is an elliptic curve $B_f$
defined over $\LL$ with CM by $\K$ and for the case that $A$ is
isogenous over $\K$ to $\Res_{\LL/\K} (B_f)$, two situations can
occur: the abelian variety $A_f$ is isogenous over $\M$ to
$\Res_{\LL/\M}(B_f)$ where  $\M$ is either  $\Q$ or $\K$. In both
cases, $\LL$ coincides with $\LL_\eta$ since $\dim (A)\leq
[\LL_\eta:\K]$ and $[\E:\Q]\, [\M :\Q]= [\LL:\Q]$.

Let $\gp$ be a prime ideal of $\K$ over $p\in\cP$. Let us denote by
$d(p)$ and $d(\gp)$ the residue degrees of $p$ and $\gp$ in $\LL$,
respectively. We have that either $d(p)=d(\gp)$ or $d(p)=2
\,d(\gp)$ depending on whether $p$ splits or it is inert in $\K$.
We know that
$$
\begin{array}{ll}
 \varepsilon (p)^{d(p)}=1 \text{ and } \psi(\gp)^{d(p)}\in\K & \text{, if $p$ splits in $\K$,}\\[3 pt]
  \varepsilon (p)^{d(p)}=1 \text{ and } \psi(\gp)^{d(p)/2}\in\K & \text{, if $p$ is inert in $\K$.}
\end{array}
$$
In both cases, one has $b_p\in\Z$. In the CM case, we introduce a
modified partition of~$\cP$:
\begin{equation}\label{s2}
\begin{array}{ll}
\cS_1=& \{ p\in\cP: d(p)=1\}\cup \{ p\in\cP: d(p)=2, \text{ $p$
inert in $\K$}\}\,,\\ [4 pt]
 \cS_2=&\{p\in\cP:
d(p)=2, \text{ $p$ splits in  $\K$}\}\,,\\[4 pt]
\cS_3=&\{ p\in\cP: d(p)\geq 3\}\,.
\end{array}
\end{equation}

Let us denote by $\cP_\K$ the set of prime ideals $\gp$ of $\K$
over all primes $p\in\cP$. Finally, let $\cS_i'$ be the subset of
$\cP_\K$ consisting on the primes over all primes in $\cS_i$.
That is,
\begin{equation}\label{s2prima}
\begin{array}{ll}
\cS_1'=& \{ \gp\in\cP_\K: d(\gp)=1\}\,,\\ [4 pt]
 \cS_{2}'=&\{\gp\in\cP_\K:
d(\gp)=2 \,,\,\,\gp\neq \overline{\gp}\}\,,\\[4 pt]
\cS_3'=&\{\gp\in\cP_\K: d(\gp)=2 \,,\,\,\gp= \overline{\gp}\}\cup
\{ \gp\in\cP_\K: d(\gp)\geq 3\}\,.
\end{array}
\end{equation}

\begin{teo}\label{teoremaCM} With the above notations, let $L_2(s):=L(f,\cS_2,s)L(f^{-},\cS_2,s)$ and
$$
L( s):=L(f,\cS_1,s)^2 L_2(s)=
\displaystyle{\left(\prod_{\gp\in\cS_1'} \frac{1}{ 1-\psi
(\gp)\Norm(\gp)^{- s}}\right)^2\prod_{\gp\in\cS_{2}'}\frac{1}{
1-\psi (\gp)^2 p^{-2 s}}}\,.
 $$
The function $L(s)^{[\LL:\K]}$ has analytic continuation to the
the right half-plane $\Re (s)>5/6$ and, moreover, one has
$$
\ord_{s=1} L(B_{f}/\LL,s)=\ord_{s=1}L(f/\LL,s)=\ord_{s=1} L(s)^{ [\LL:\K]} \,,
$$
and, thus,
$$\ord_{s=1} L(A_f/\LL,s)=\displaystyle{[\E:\Q]\ord_{s=1} L(s)^{ [\LL:\K]}\,.}
$$
Moreover, if $A_f$ is isogenous over $\M$ to
$\Res_{\LL/\M} (B_f)$ for $\M\subseteq \K$, then
$$
\ord_{s=1} L(A_f/\Q,s)=\displaystyle{\frac{1}{[\M:\Q]}\ord_{s=1}
L(s)^{[\LL:\K]}}\,.
$$
\end{teo}
\noindent {\bf Proof.} Observe that $L(A_f/\K,s)=L(A_f/\Q,s)^2$.
Now, we can factorize $L^*(f/\LL,s)=\prod_{i=1}^3 G_i(f,s)$, where
$$
\begin{array}{ll}
G_1(f,s):=&\displaystyle{L(f,\cS_1,s)^{[\LL:\Q]}},\\[4 pt]
G_2(f,s):=&\displaystyle{\left(\prod_{p\in\cS_2}\frac{1}{ 1-b_p p^{-
2\,s}+
 p^{2 (-2 s+1)}}\right)^{[\LL:\Q]/2}},\\[4 pt]
G_3(f,s):=&\displaystyle{\left(\prod_{p\in \cS_3}\frac{1}{ 1-b_p
p^{- d(p)\,s}+ p^{d(p) (-2
s+1)}}\right)^{[\LL:\Q]/d(p)}}.
\end{array}
$$
By taking into account that
$L^*(A_f/\LL,s)=L^*(f/\LL,s)^{[\E:\Q]}$ and that if
$\Res_{\LL/\K}(B_f)$ is isogenous to $A_f$ over $\K$  then
$\Res_{\LL/\Q}(B_f)$ is isogenous to $A_f^2$ over $\Q$, the
statement is obtained by using similar arguments as in the proof
of Theorem~\ref{teoremanoCM}. \hfill$\Box$

\vskip 0.5truecm

Notice that the same arguments used in Remark~\ref{rem2} can be
applied to obtain the corresponding analogues for the CM case. As
a by-product, we also obtain the following application to the
elliptic curves studied by B.\,Gross in~\cite{gross}.

\begin{prop}\label{gross}
Let $p\equiv 3 \pmod 4$ be a rational prime $>3$. Let $A(p)$ be
the Gross's elliptic curve with CM by the ring of integers of
$\K=\Q(\sqrt{-p})$. Let $f$ be a normalized newform  in
$S_2(\Gamma_0(p^2))$ such that $A(p)$ is a quotient of~$A_f$.
Then,
$$\ord_{s=1}L(f,s)=\frac{1}{[\HH:\Q]}\ord_{s=1} L(f,\cS_1,s)^{[\HH:\Q]}\,,$$
where $\HH$ denotes the Hilbert class field of $\K$.
\end{prop}

 \noindent
{\bf Proof.} In this case, $\LL=\HH$ and the class
number of $\K$, say $h$, coincides with $[\LL:\K]$ and $[\E:\Q]$.
Moreover, $\cS_2=\emptyset$ since  $h$ is odd and, thus, the
function $L(s)$ in the above theorem is $L(f,\cS_1,s)^2$.
We know  that, for all $\sigma\in\Gal(\overline{\Q}/\Q)$, the value $\ord_{s=1} L({}^\sigma f,s)$
is either $0$ or $1$ depending on whether $p\equiv 7\pmod 8$ or $p\equiv 3\pmod 8$
(cf. \cite{Mo-Ro} and \cite{yang2000}).  To prove the statement,
we use that  $\Res_{\LL/\K}(A(p))$ is isogenous to $A_f$ over $\K$
and that $\ord_{s=1} L(f,s)$ is invariant under Galois
conjugations of $f$. Then, applying  Theorem~\ref{teoremaCM} we
obtain
$$[\E:\Q]\ord_{s=1}L(f,s)= \ord_{s=1} L(A_f/\Q,s)=  \frac{1}{2} \ord_{s=1}L(f,\cS_1,s)^{[\LL:\Q]} \,.$$
The statement is now an immediate consequence. \hfill $\Box$

\section{Distributions of Frobenius traces and their convergences in average}

With the aim of enlightening that the role of primes $p$ having residue
degree $2$ in $\LL$ is minor, in this section we examine the convergence of
the product
 $$
 \prod_{p\in\cS_2}\frac{1}{ 1-b_p p^{- 2}+
 p^{-2}}=
\prod_{p\in\cS_2}\left( 1- \frac{b_p-1}{p^2}\right)^{-1}
\,,
 $$
for any newform $f$ (with or without CM). To this end, first we
generalize some results of T.\, Mitsui for the CM case and we
prove that one can spurn these primes. As for the non-CM case,
we adapt the proof of the result obtained by T. Barnet-Lamb, D. Geraghty,
M. Harris and R. Taylor in \cite{BGHT} about Sato-Tate distributions to extend it
to the case when the set of primes is restricted to an arithmetic
progression and, then, we prove an estimation for the rate of
convergence of the mathematical expectations.

In the sequel we will use the following notation. For any subset
$\cS$ of $\cP$ and a real number $t$, we put $\cS(t):=\{
p\in\cS:p\leq t\}$ and $\delta_t$ denotes the Dirac measure at
$t$. Before going further, we need the following two lemmas.

 \begin{lema}\label{convergence}
Let $\cS=\{p_1<\cdots< p_n < \cdots\}$ be an infinite sequence of
rational primes and let $\{c_n\}_{n>0}$  be a  sequence of real
numbers such that
 $\mid c_n\mid \leq k p_n$ for some  $k>0$  and all $n>0$. Then, the product  $
 \prod_{n>0}( 1-c_n p_n^{- 2})^{-1}
 $ converges  if and only if $\sum_{n>0}c_n/p_n^2$ converges.
 \end{lema}

\noindent {\bf Proof.} For every positive real $t$ we set
 $$
 \mathfrak S(t)=\prod_{p_n\in \cS(t)}\frac{1}{ 1-c_n p_n^{- 2}} \,.
 $$
The convergence of $\prod_{n>0}( 1-c_n p_n^{- 2})^{-1}$
is equivalent to the existence of $\lim_{t\to
+\infty}\log\mathfrak S(t)$. Now,
$$
\log\mathfrak S(t)=-\sum_{p_n\in \cS(t)}\log\left(1-\frac{c_n}{p_n^2}\right)=
-\sum_{p_n\in \cS(t)}\left(\log\left(1-\frac{c_n}{p_n^2}\right)+\frac{c_n}{p_n^2}\right)+\sum_{p_n\in \cS(t)}\frac{c_n}{p_n^2}\,,
 $$
 and the convergence of the first term follows from
 $$
 \sum_{p_n \in \cS(t)}\left| \log\left(1-\frac{c_n}{p_n^2}\right)+\frac{c_n}{p_n^2}\right|
 < \sum_{p_n \in \cS(t)}  c_n^2/p_n^4 < \infty
 $$
 since  $c_n^2/p_n^4 <k^2/p_n^2$.\hfill $\Box$

 \vskip 0.2cm

 \begin{lema}\label{expectation}
Let  $\cS=\{ p_1< \cdots< p_n <\cdots \}$ be a sequence  of primes
of positive density and let $\{ c_n\}$ be a sequence of real
numbers. Consider the  sequence of  the probability measures
$\nu_n:=1/n \sum_{i=1}^n \delta_{c_{i} / p_i}$ and assume that the
sequence of their mathematical expectations, $E(\nu_n)$, converges
to a real number $\ell$.
  Then, we have that
  \begin{itemize}
 \item[(i)] if $\ell\neq 0$, then the series  $\sum_{n>0} c_n/p_n^2$ does not converge,

\item[(ii)] if $\ell=0$ and the real function $\cK(t)= \frac{1}{|\cS(t)|}\sum_{p_n\in\cS( t)}c_n/p_n$
defined on $[p_1,+\infty)$ satisfies the condition
$\int_{p_1}^{+\infty} \frac{\mid \cK (t)\mid }{t\, \log
t}\,dt<+\infty$, then the series $\sum_{n>0} c_n/p_n^2$ is
convergent.
  \end{itemize}
\end{lema}

\noindent
{\bf Proof.} By hypothesis, we have
$$
\lim _{n\to +\infty} E(\nu_n)=\lim _{t\to+\infty} \frac{1}{|\cS(t)|}\sum_{p_n\in\cS(t)} \frac{c_n}{p_n}=\ell.
$$
 Let $\cK(t)$ be the function defined by the relation
$$
\sum_{p_n\in\cS( t)}\frac{c_n}{p_n}=|\cS(t)| \ell+|\cS (t)| \cK(t)\,.
$$
By partial summation  we get
$$
\displaystyle{\sum_{p_n\in \cS(t)}\frac{c_n}{p_n^2}=\frac{1}{t}\sum_{p_n\in\cS( t)}\frac{c_n}{p_n}+\int_{p_1}^t \frac{1}{x^2}
\sum_{p_n\in\cS( x)}\frac{c_n}{p_n}\,dx.}
$$
Observe that $|\cS(t)|=c\, t/\log t(1+ o(1))$ and $p_n=  c\, n\,
\log n(1+ o(1))$, where $c>0$ is the density of $\cS$. On one
hand, we have that
$$
\lim_{t\to+\infty} \frac{1}{t}\sum _{p_n\in\cS(t)}\frac{c_n}{p_n}=\lim_{t\to+\infty} \left(\frac{1}{|\cS(t)|}\sum _{ p_n\in\cS(t)}\frac{c_n}{p_n}\right)\frac{|\cS(t)|}{t}= \ell\,\lim_{t\to+\infty}\frac{|\cS(t)|}{t}=0\,.
$$
On the other hand, one has that
\begin{equation}\label{eqk}
\displaystyle{\int_{p_1}^t \frac{1}{x^2}
\sum_{p_n\in\cS(x)}\frac{c_n}{p_n}\,dx= \ell\,\int_{p_1}^t \frac{|\cS (x)|}{x^2}\, dx+\int_{p_1}^t \frac{|\cS (x)|}{x^2} \cK (x)\, dx\,.}
\end{equation}
Since $\lim_{t\to +\infty}\cK (t)=0$ and $\lim_{t\to +\infty}\int_{p_1}^t|\cS (x)|/x^2\,dx=+\infty$,
the convergence of $\sum_{n>0} c_n/p_n^2$ implies that $\ell=0$.

Assume now that $\ell=0$. By the hypothesis in (ii), the last integral in (\ref{eqk})
is absolutely convergent and, thus, the series $\sum_{i>0}c_n/p_n^2$ converges.
\hfill $\Box$

\begin{rem}
Recall that if a sequence of probability measures $\nu_n$ with
support on a compact real subset converges in law to a continuous
probability measure $\nu$, then the mathematical expectation
$E(\nu)$ exists and the sequence of the mathematical expectations
$E(\nu_n)$ converges to~$E(\nu)$. Moreover, if a real function
$\cK$  provides a uniform estimation for this convergence, namely
$\nu_n( (-\infty, a])-\int_{-\infty}^a d\nu =O( \cK (n))$ for all
real~$a$,  then $E(\nu_n)-E(\nu)=O( \cK (n))$. We will apply the
above lemma in this sense, proving that the sequence $\nu_n=\sum_{p\in\cS_2(n)}
\delta_{b_{p}/(2 p)}$ converges in law to a continuous
probability measure with support on $[-1,1]$.
\end{rem}

Next, we treat in separate subsections the CM  and non-CM cases.

\subsection{Uniform distribution of arguments for the CM case}
In this subsection, it is assumed that $f$ has CM by an imaginary
quadratic field  $\K$ and we keep the notation as in
Section~\ref{CM}.

For the particular case that $A_f$ is an elliptic curve, then
$\LL=\K$ and  it is well-known that the values $a_p/(2\sqrt{p})\in
\R $ for the set of  primes $p\in\cP$ which split in $\K$, say
$\cP_1$, are equidistributed in $[1, 1]$ with respect to the
measure  $d\mu = 1/\pi(1-x^2)^{-1/2}\,dx $; i.e., for all
$a\in[-1,1]$ we have that
\begin{equation}\label{uniform}
\displaystyle{\lim_{t\to +\infty} \frac{|\{p\in \cP_1(t):
\frac{a_p}{2\sqrt p}\leq a \}|}{
|\cP_1(t)|}= \frac{1}{\pi}\int_ {-1}^a \frac{1}{\sqrt{1-x^2}}\, dx}\,.
\end{equation}
In other words, the sequence of measures $\frac{1}{n}\sum_{i=1}^n
\delta_{ a_{p_i}/(2 \sqrt{p_i})}$ converges in law to $\mu$, where
$\cP_1=\{p_1<\cdots<p_n<\cdots\}$. If we take $\theta_p\in (0,
\pi)$ such that $a_p/(2\sqrt p)=\cos \theta_p$, then the above
condition amounts to saying that the sequence of measures
$\frac{1}{n}\sum_{i=1}^n \delta_{ \theta_{p_i}}$ converges in law
to the uniform measure on $[0, \pi]$.

In order to remove the set of primes $\cS_2$ from Theorem~\ref{teoremaCM}, we shall
generalize~(\ref{uniform}), in Theorem~\ref{jorge}, to CM modular abelian
varieties $A_f$ of higher dimension, for several subsets of primes which
split completely in $\K$, and also including an estimate for the error term.
The next result is the main tool to achieve this goal.

\begin{teo}\label{quadrat} Let $\psi'$ be any Hecke character of $\K$
mod $\gn$ (not necessarily primitive), and let
$\theta\in[0,2\pi]$. Then, for any ideal class $C$ in
$I(\gn)/P_1(\gn)$ one has
\begin{equation}\label{part1}
| \{ \gp\in C\cap\cP_\K: \Norm (\gp)\leq t, \arg \psi'(\gp)\leq
\theta\}| =\frac{\theta}{2\pi\,[\K_{\gn}:\K] }\int_2^{t}
\frac{dx}{\log x}+ O\left({t} \,e^{-c\sqrt{\log t}}\right)\,,
\end{equation}
where $c$ is some positive real constant. In particular,
\begin{equation}\label{part2}
\displaystyle{ \frac{| \{ \gp\in C\cap\cP_\K: \Norm (\gp)\leq t,
\arg \psi'(\gp)\leq \theta\}|} {| \{ \gp\in C \cap\cP_\K: \Norm
(\gp)\leq t \}|}-\frac{\theta}{2\pi}=O\left( e^{-K\sqrt{\log
t}}\right)}
\end{equation}
for some positive real constant $K$.
\end{teo}

\noindent {\bf Proof.} The proof is based on a theorem of
T.\,Mitsui in \cite{mitsui} and uses similar arguments. Let us
note that it is enough to restrict to primitive characters, since
the only difference will be in a finite number of primes dividing
the conductor. Now, first, let  $\theta\in [0,2\pi)$, $t>2$ and
$\Delta=e^{-u\sqrt{\log t}}$ for some $u$ that will be chosen
later. We will denote $C(t)=\{\ga\in C\,:\,\Norm(\ga)\le t\}$.
Consider a  $2\pi$-periodic function $g\in {\mathcal C}^2(\mathbb
R)$ such that $g(x)=1$ for $x\in[0,\theta]$, $0\le g(x)\le1$ for
$x\in[-\Delta,0]\cup[\theta,\theta+\Delta]$ and $g(x)=0$ for any
other $x\in[-\pi+\frac{\theta}2,\pi+\frac{\theta}2]$. Assume also that
$||g'||_2^2\ll 1/\Delta$ and let $g(x)=\sum_nc_ne^{inx}$ denote its
Fourier expansion. For any given ideal $\ga$, let $\theta_{\ga}$
be the argument of $\psi'(\ga)$; in other words,
$\psi'(\ga)=\sqrt{\Norm(\ga)}e^{i\theta_{\ga}}$. Hence, summing
over prime ideals in $C(t)$, we get
\begin{equation}\label{fourier}
\#\{\gp\in C(t)\,:\, 0\le \theta_{\gp}\le \theta\}\le\sum_{\gp\in
C(t)}g(\theta_{\gp}).
\end{equation}
First, we  will  proceed to obtain the right asymptotics for the
sum. This will provide an upper bound for the left hand side of
(\ref{fourier}). An analogous argument will give us the lower
bound and, hence, the theorem. By expanding into its Fourier
series we get
\begin{eqnarray*}
\sum_{\gp\in C(t)}g(\theta_{\gp})&=&\sum_{\gp\in
C(t)}\sum_{n}c_ne^{i n\theta_{\gp}}=
\frac{1}{2\pi}\int_0^{2\pi}g(x)dx\sum_{\gp\in C(t)}1+\sum_{n\ne
0}c_n\sum_{\gp\in C(t)}e^{i n\theta_{\gp}} \,.
\end{eqnarray*}
Let us denote the two terms by
$$
T_1 = \frac{1}{2\pi}\int_0^{2\pi}g(x)dx\sum_{\gp\in C(t)}1
\,,\qquad T_2 = \sum_{n\ne 0}c_n\sum_{\gp\in C(t)}e^{i
n\theta_{\gp}}\,.
$$
Noting that $\int_0^{2\pi}g(x)dx=\theta+O(\Delta)$, Theorem in
Section~1 of~\cite{mitsui} gives us for the main term
$$
T_1=\frac{\theta}{2\pi[\K_{\gn}:\K]} \,\li{t}+O(\Delta\, \li{t}),
$$
where, as usual, $\li{t}=\int_2^t\frac1{\log x}dx$ is the
logarithmic integral. Observe that, for any choice of $u$ in
$\Delta$, this error term is of the form stated in (\ref{part1})
and then, to prove the first part of the theorem, we have to get
the correct upper bound for the error term $T_2$. Let $\xi:
I(\gn)\to\C^*$ be the group homomorphism defined by
$$
\xi(\ga):=\frac{\psi'(\ga)}{\sqrt{\Norm(\ga)}}=e^{ i\theta_ {\ga}}.
$$
By definition, for any integer $n\ne0$, the map $\xi^n$ is a
non-real Grossencharacter of the field $\K$ of modulus $\gn$ and
frequency $n$ as in Section 3.8 of \cite{iwaniec}. In particular,
$$
L(\xi^n,s)=\sum_{(\ga,\gn)=1}\frac{e^{ in\theta_
{\ga}}}{\Norm(\ga)^{s}}
$$
is a Hecke $L$-function. With the notation  as in Chapter $5$ of \cite{iwaniec}, this $L$-function
has parameters $d=2$, $\kappa_1=|n|/2$, $\kappa_2=|n|/2+1$ and
$\mathfrak q(\xi^n)\ll |n|^2$.
Hence, by noting that
$$
\sum_{\gp\in C(t)}e^{ i n\theta_{\gp}}=\sum_{m\le t}\frac{\Lambda_{\xi^n}(m)}{\log m}+O\left(\sqrt t\log t\right),
$$
where here $\Lambda_{\xi^n}(m)= (e^{i kn\theta_{\gp}}+e^{i
kn\theta_{\bar{\gp}}})\log p$ for any $m=p^k$ with $\Norm(\gp)=p$, and
$\Lambda_{\xi^n}(m)=0$ otherwise, and that the constant in the error
term is independent of $n$, we can apply Theorem~5.13 of
\cite{iwaniec} and partial summation to obtain
\begin{equation}\label{eq:t2}
T_2\ll t\sum_{n\ne 0}|c_n|\exp\left(-c\frac{\log t}{\sqrt {\log t}+6\log |n|}\right).
\end{equation}
Now,  the Cauchy-Schwarz inequality yields
\begin{equation}\label{eq:coeff}
\left(\sum_{n\in S}|c_n|\right)^2\le\sum_{n\in S}n^2|c_n|^2\sum_{n\in S}\frac1{n^2}\ll||g'||_2^2\sum_{n\in S}\frac1{n^2}\ll \frac1\Delta\sum_{n\in S}\frac1{n^2}
\end{equation}
for any set $S\subset \mathbb Z$. Hence, letting $R=e^{v\sqrt{\log
t}}$, for some $v$ to be chosen later,  and splitting the sum in
the right hand side of (\ref{eq:t2}) into two different sums
$T_{2,1}$, $T_{2,2}$, depending on whether $|n|<R$ or not
respectively, we get on one hand
$$
T_{2,1}\le te^{-\frac{c}{6v+1}\sqrt{\log t}}\sum_{|n|<R}|c_n|\ll t\frac1{\Delta^{1/2}}e^{-\frac{c}{6v+1}\sqrt{\log t}}=t e^{(\frac u2-\frac{c}{6v+1})\sqrt{\log t}}.
$$
On the other hand, using the trivial bound $\exp\left(-c\frac{\log
t}{\sqrt {\log t}+6\log |n|}\right)<1$ and (\ref{eq:coeff}),  we
get for $|n|\geq R$,
$$
T_{2,2}\ll t\frac{ 1}{(\Delta R)^{1/2}}=e^{(u-v)/2\sqrt {\log t}}.
$$
From here  it is easy to  get the upper bound
$$
T_2=O\left(te^{-c'\sqrt {\log t}}\right),
$$
for some $c'>0$,  for example by choosing $(6v+1)v=c/2$ and
$u=v/2$.

In order to get the lower bound, we just have to consider a
$2\pi$-periodic function $\hat g\in {\mathcal C}^2(\mathbb R)$
such that $\hat g(x)=1$ in the interval $[\Delta,\theta-\Delta]$,
$0\le \hat g(x)\le 1$ for $x\in
[0,\Delta]\cup[\theta-\Delta,\theta]$ and $\hat g(x)=0$ for any
other $x\in [-\pi+\frac{\theta}{2},\pi+\frac{\theta}{2}]$. We can assume
that~$||\hat g'||^2\ll \frac1{\Delta}$, for the same selection of
$\Delta$ as above. Now,
$$
\sum_{\gp\in C(t)}\hat g(\theta_{\gp})\le\#\{\gp\in C(t)\,:\, 0\le \theta_{\gp}\le \theta\},
$$
and the same arguments as above produce
$$
\sum_{\gp}\hat g(\theta_{\gp})=\frac{\theta}{2\pi [\K_{\gn}:\K]} \,\li{t}+O\left(te^{-c\sqrt{\log t}}\right)\,,
$$
for some constant $c$, which ends the proof of (\ref{part1}) for
$\theta\in [0,2\pi)$. The case $\theta=2\pi$ is given by Theorem
in Section~1 of \cite{mitsui}. Finally, observe that~(\ref{part2})
is an immediate consequence of~(\ref{part1}).~\hfill~$\Box$


\begin{rem}
We observe that for the case when $C$ is a class of principal
ideals, the above theorem agrees with Theorem in Section~4 of
\cite{mitsui} adapted to imaginary quadratic fields. Indeed, one
has $C=P_\delta(\gn)$ for some $\delta\in (\cO/\gn)^*$ and let
$\eta'$ be the corresponding character mod $\gn$ attached to
$\psi'$. For any $\alpha\in\cO$ such that $\alpha\equiv
\delta\pmod{\gn}$, we have $\psi'( \alpha
\cO)=\alpha\,\eta'(\delta)$ and, thus, $\arg \psi'(\alpha
\cO)\equiv \arg \alpha+ \arg \eta'(\delta) \pmod{2\pi}$.
\end{rem}

\vskip 0.2 truecm

As a consequence, we obtain the following result concerning the CM
newform~$f$. Recall that for a prime $p\in\cP$ (resp.
$\gp\in\cP_{\K}$), $d(p)$ (resp. $d(\gp)$) denotes its residue
degree in $\LL$, and let $b_p$ be as in the beginning of
Section~2.
\begin{teo}\label{jorge} Let $n>0$ be an integer such that the set
$$
\cP_n:=\{ p\in\cP:  d(p)=n\,,\text{ $p$ splits in $\K$ }\}
$$
is non-empty. For every $a\in [-1,1]$, one has
$$
\displaystyle{ \frac{|\{p\in \cP_n(t):
\frac{b_p}{2  p^{n/2}}\leq a \}|}{
|\cP_n(t)|}- \int_ {-1}^a d\mu =O \left(  e^{-K\sqrt{\log t}}\right)}\,,
$$
for some positive constant $K$.

\end{teo}
\noindent {\bf Proof.} Since the conductor of  $\varepsilon$ and
$\Norm (\gm)$ divide the level $N$ of the newform $f$ and $\LL/\K$
is an abelian extension, the field $\LL$ is contained in the ray
class field of $\K$ mod $\gn$,  where $\gn=N\cO$.
Consider the sets
$$\cP_{n,\K}:=\{\gp\in \cP_\K\colon d(\gp)=n\}\quad\text{and}\quad \cP_{n,\K}':=\{\gp\in \cP_{n,\K}\colon\gp\neq \overline{\gp}\}\,.$$
Note that $\cP_{n,\K}'$ is the subset of  $\cP_{\K}$ consisting on
the primes over all  primes in $\cP_n$. It is clear that if
$\gp\in \cP_{n,\K}$, then any prime  of $\K$ in the class of $\gp$
in $I(\gn)/P_1(\gn)$ also lies in $\cP_{n,\K}$. Hence, there
exists a subset $S$ of $I(\gn)$ which is a union of certain classes of
$I(\gn)/P_1(\gn)$ and such that $S\cap\cP_\K=\cP_{n,\K}$.
 For every class of $S$, we apply~(\ref{part1}) of Theorem~\ref{quadrat},
replacing $\psi'$ with $\psi$ as Hecke character mod  $\gn$ and
then summing over the set of classes in $S$, and we obtain
$$
\displaystyle{
\frac{| \{ \gp\in  \cP_{n,\K} : \Norm (\gp)\leq t, \arg \psi(\gp)\leq \theta\}|}
{| \{ \gp\in   \cP_{n,\K} : \Norm (\gp)\leq t \}|}-\frac{\theta}{2\pi}=O\left( e^{-K\sqrt{\log t}}\right)}\,.
$$
The set $ \cP_{n,\K}$ is the disjoint union of $\cP_{n,\K}'$ and
the set of the ideals generated by inert primes $p$ in $\K$ such
that $p\cO\in \cP_{n,\K}$. Since the density of this last set in
$\cP_{n,\K}$ is zero,  we get
 $$
\displaystyle{
\frac{| \{ \gp\in \cP_{n,\K}': \Norm (\gp)\leq t, \arg \psi(\gp)\leq \theta\}|}
{| \{ \gp\in  \cP_{n,\K}': \Norm (\gp)\leq t \}|}-\frac{\theta}{2\pi}=O\left( e^{-K\sqrt{\log t}}\right)}\,,
$$
and, thus, also
 $$
\displaystyle{
\frac{| \{ \gp\in \cP_{n,\K}': \Norm (\gp)\leq t, \arg \psi(\gp^n)\leq \theta\}|}
{| \{ \gp\in  \cP_{n,\K}': \Norm (\gp)\leq t \}|}-\frac{\theta}{2\pi}=O\left( e^{-K\sqrt{\log t}}\right)}\,.
$$
 For $p\in\cP_n$, the two primes $\gp,\overline{\gp}\in \cP_{n,\K}'$ over $p$  satisfy
$\arg \psi (\gp^n)+\arg \psi (\overline{\gp}^n)=2\pi$. Let $\gp$
be the unique prime such that $\arg \psi (\gp^n)\in (0,\pi)$.
Hence, $b_p=2\, p^{n/2} \, \cos (\theta_p)$ with $\theta_{p}=\arg
\psi (\gp^n)$ and it follows the statement. \hfill $\Box$

\begin{rem}
Note that  the above theorem for $n=1$ generalizes and improves
(\ref{uniform}) for higher dimensional CM modular abelian varieties
$A_f$ and for $n=2$ states
$$
\displaystyle{ \frac{|\{p\in \cS_2(t):
\frac{b_p}{2  p}\leq a \}|}{
|\cS_2(t)|}- \int_ {-1}^a d\mu =O \left(  e^{-K\sqrt{\log t}}\right)}\,,
$$
for some positive constant $K$.
\end{rem}

Next, we present the main theorem for the CM case.

 \begin{teo}\label{mainCM}
Assume that $f$ has CM and let $E/\LL $ be an elliptic curve such
that $A_f$ is isogenous over $\LL$ to $E^{\dim A_f}$. Then, one
has
 $$
 \ord_{s=1} L(E/\LL,s)=\ord_{ s=1} L(f/\LL,s)=\ord_{ s=1} L(f,\cS_1,s)^{[\LL:\Q]}\,.$$
 Moreover, if $A_f$ is isogenous over $\M$ to
$\Res_{\LL/\M} (E)$ for $\M\subseteq \K$, then
$$
\ord_{s=1} L(A_f/\Q,s)=\displaystyle{\frac{1}{[\M:\Q]}\ord_{s=1}
L(f,\cS_1,s)^{[\LL:\Q]}}.
$$
 \end{teo}

\noindent
 {\bf Proof. }
Since $E(\mu)=0$, the statement follows from
Lemmas~\ref{convergence},~\ref{expectation}, Theorem~\ref{jorge}
and Theorem~\ref{teoremaCM}.\hfill $\Box$

\begin{coro}  With the same notation as in Theorem~\ref{mainCM},
assume that $A_f$ is isogenous over $\M$ to $\Res_{\LL/\M} (E)$
for $\M\subseteq \K$. If $\ord_{ s=1} L(f,s)$ is invariant under
Galois conjugations of $f$, then we have
$$
\ord_{ s=1} L(f,s)=\frac{1}{[\LL:\Q]} \ord_ {s=1} L(f,\cS_1,s)^{[\LL:\Q]}\,.
$$
\end{coro}

This corollary is obtained as a consequence of the above
discussion and it generalizes our previous result concerning
Gross's elliptic curves in Proposition~\ref{gross}.

\begin{rem} In the above corollary, it is assumed that
$\ord_{s=1}L(f,s)$ is invariant under Galois conjugations acting
on $f$. Although this is expected to be always true, only a few
results are known in this direction (see Corollary 1.3 of
Gross-Zagier in \cite{grosszagier}).
\end{rem}

\begin{rem}\label{rank}
For any number field $\M$, let $r_\M=\rank(A_f(\M))$ be the rank
of the Mordell-Weil group over $\M$. From the inclusion
$\E\hookrightarrow \End (A_f(\M)\otimes \Q)$ when $r_\M>0$, it
follows that $A_f(\M)\otimes \Q$ is also an $\E$-vector space and,
thus, $[\E:\Q]$ divides  $r_\M$.  According to the  Birch and
Swinnerton-Dyer conjectures, one expects   that $[\E:\Q]$ divides
$\ord_{s=1} L(A_f/\Q, s)$. In the CM case,  if $A_f$ is isogenous
over $\M$ to $\Res_{\LL/\M} (E)$ for $\M\subseteq\K$, then this
divisibility condition implies that the rational number
$$\displaystyle{\frac{1}{[\LL:\Q]}\ord_{s=1}
L(f,\cS_1,s)^{[\LL:\Q]}}$$ must be an integer.
\end{rem}

\subsection{Sato-Tate distributions for the non-CM case}

In this subsection, we assume that $f$ does not have CM. First, we include
a result about  Sato-Tate distributions in arithmetic progressions.
Let $\mu$ be the so-called Sato-Tate measure; i.e., the measure
with support on $[-1,1]$ such that $d\mu=2/\pi\sqrt{1-x^2}dx$. The next result generalizes
part~3 of Theorem~B of~\cite{BGHT}.

\begin{teo}
Let $M$ be any multiple of the conductor of $\varepsilon$. Let $\zeta$ be a root of unity such that $\zeta^2=\varepsilon(m)$ for some $m\in (\Z/M\,\Z)^*$. Then, for all $a\in[-1,1]$ we have
$$
\lim_{t\to +\infty}\frac{|\{ p\in\cP (t):  \frac{a_p}{2\sqrt p \,\,\zeta}\leq a\,,\,\, p\equiv m\pmod M|}
{|\{ p\in\cP(t):   p\equiv m\pmod M|}=\int_{-1}^a d\mu\,.
$$
\end{teo}

{\bf Proof.} For a prime $p\in\cP$, let us denote by $\alpha_p$ and $\beta_p$ the roots of the polynomial $x^2-a_p \,x+p\,\varepsilon(p)$. We know that the statement is equivalent to prove that for all integer $n>0$ one has
$$
\sum_{\substack{p\in\cP( x) \\ p\equiv m\,(\bmod M)}}\frac{ \sum _{i=0}^n \alpha_p^{n-i}\beta_p^i}{p^{n/2}\zeta^n}= o\left(\frac{x}{\log x}\right)\,.
$$
Of course, $\zeta^n$ can be omitted in the above condition. Let $\chi$ be a Dirichlet character mod $M$.
Consider the partial $L$-function
$$
L((\operatorname{Symm}^n\, f)\otimes \chi,s)=\prod_{p\in\cP} \prod_{i=0}^n (1-\alpha_p^{n-i}\beta_p ^i\, \chi(p)/p^s)^{-1},
$$
which converges absolutely on $\Re( s)>n/2+1$. By part 2 of Theorem B in \cite{BGHT}, we know that this function has meromorphic continuation to the whole complex plane and is holomorphic and non-zero in $\Re(s)\geq n/2+1$. Therefore, its logarithmic derivative
$$
H_{\chi,n}(s):=\frac{d}{ds} \log\left(L((\operatorname{Symm}^n\, f)\otimes \chi,s)\right)=\frac{\frac{d}{ds}L((\operatorname{Symm}^n\, f)\otimes \chi,s)}{L((\operatorname{Symm}^n\, f)\otimes \chi,s)}
$$
is holomorphic for $\Re( s)\geq n/2+1$. Hence,
$$
\begin{array}{ll}
H_{\chi,n}(s)&=-\displaystyle{\sum_{p\in\cP}\sum_{k\geq 1} \frac{\sum_{i=0}^n(\alpha_p^{n-i}\beta_p ^i\, \chi(p))^k\,\log p}{p^{k\, s}}}\\ &=-\displaystyle{\sum _{p\in\cP}\frac{\sum_{i=0}^n \alpha_p^{n-i}\beta_p ^i\, \chi(p) \log p}{p^s}+ G(s),}
\end{array}
$$
where $G(s)$ is holomorphic for $\Re (s)>n/2+1/2$. Summing over all Dirichlet characters $\chi$ mod $M$,  we get:
$$\displaystyle{\sum_{\chi} \sum _{p\in\cP}\overline{\chi(m)}\frac{\sum_{i=0}^n \alpha_p^{n-i}\beta_p ^i\, \chi(p) \log p}{p^s}}=-\varphi(M)\displaystyle{\sum_{\substack{p\in\cP( x) \\ p\equiv m\,(\bmod M)}}\frac{\sum_{i=0}^n \alpha_p^{n-i}\beta_p ^i\, \log p}{p^s}}\,,
$$
where $\varphi$ stands for Euler's function. Since the function
$$
\displaystyle{\sum_{\substack{p\in\cP( x) \\ p\equiv m\,(\bmod M)}}\frac{\sum_{i=0}^n \alpha_p^{n-i}\beta_p ^i\, \log p}{p^{n/2}p^s}}
$$
is holomorphic for $\Re(s)\geq 1$, by the Wiener-Ikehara tauberian theorem (cf.
 Theorem~1 in Chapter  XV $\S$3  of \cite{lang}),  we obtain that
$$\displaystyle{\sum_{\substack{p\in\cP( x) \\ p\equiv m\,(\bmod M)}}\frac{\sum_{i=0}^n \alpha_p^{n-i}\beta_p ^i\, \log p}{p^{n/2}}}=o(x).
$$
Now, applying Abel summation, we get
$$\sum_{\substack{p\in\cP( x) \\ p\equiv m\,(\bmod M)}}\frac{\sum_{i=0}^n \alpha_p^{n-i}\beta_p ^i}{p^{n/2}}=o\left(\frac{x}{\log x}\right),
$$
as we wanted to  prove. \hfill$\Box$

\begin{rem}
Note that if the Nebentypus $\varepsilon$ of the newform $f$ is
trivial, then Sato-Tate distribution applies to arbitrary
arithmetic progressions.
\end{rem}

Let $\cS_2$ be the set of primes defined in (\ref{s1}). Assume that
$\cS_2$ is non-empty, and consider the following subsets of
$\cS_2$:
 $$
\begin{array}{lll}
\cS_2^+&=&\{ p\in\cS_2: \varepsilon(p)=\phantom{-}1\}\,,\\
\cS_2^-&=&\{ p\in\cS_2: \varepsilon(p)=-1\}\,.
\end{array}
$$
Clearly, $\cS_2$
 is the disjoint union of $\cS_2^+$ and $ \cS_2^-$.
Since each of them is the set of all primes in certain congruence
classes mod $N$, as a consequence of the above theorem we obtain the following result.

\vskip 0.2 cm

\begin{coro}\label{una}
 For all $a\in[-1,1]$ we have that
\begin{itemize}
\item [(i)] if $\cS_2^+$ is non-empty, then
$$
\displaystyle{\lim_{t\to +\infty} \frac{|\{p\in \cS_2^+(t):
\frac{a_p}{2\sqrt p}\leq a \}|} {|\cS_2^+(t)| }=\int_ {-1}^a d\mu
}\,,
$$
\item[(ii)] if $\cS_2^-$ is non-empty, then
$$
\displaystyle{\lim_{t\to +\infty} \frac{|\{p\in \cS_2^-(t):
\frac{a_p}{2\, i\,\sqrt {p}}\leq a\}|}{ |\cS_2^-(t)|}= \int_
{-1}^a d\mu}\,.
$$
\end{itemize}
\end{coro}

We introduce the  measures $\mu^+$ and $\mu^-$ with support on
$[-1,1]$ and defined as follows:
$$\begin{array}{lr}\mu^+((-\infty,x\,]):=&\mu ([-\sqrt{(x+1)/2},\sqrt{(x+1)/2}\,\,])\,,\\[ 6pt]
\mu^-((-\infty,x\,]):=& 1- \mu
([-\sqrt{(1-x)/2},\sqrt{(1-x)/2}\,\,])\,;
\end{array}$$
that is,
$$
\displaystyle{d\mu^+= \frac{1}{\pi}\sqrt{\frac{1-x}{1+x}}\,
dx\quad \text{and}\quad
d\mu^-=\frac{1}{\pi}\sqrt{\frac{1+x}{1-x}}\, dx\,.}
$$

Recall that for $p\in\cS_2$, we have $b_p=a_p^2-2\,p\varepsilon(p)$. As a by-product of the above corollary, it follows the next result.

\begin{coro}\label{corouna}
Assume $\cS_2$ is non-empty. For all  $a\in [-1,1]$  we have that
\begin{itemize}
\item [(i)] if $\cS_2^+$ is non-empty, then
$$
\displaystyle{\lim_{t\to +\infty} \frac{|\{p\in \cS_2^+(t):
\frac{b_p}{2 p}\leq a\}|} {|\cS_2^+(t)|}=\int_ {-1}^a d\mu^+ }\,,
$$
\item[(ii)]  if $\cS_2^-$ is non-empty, then
$$
\displaystyle{\lim_{t\to +\infty} \frac{|\{p\in \cS_2^-(t):
\frac{b_p}{2 p}\leq a\}|} {|\cS_2^-(t)|}=\int_ {-1}^a d\mu^- }\,.
$$
\end{itemize}
\end{coro}

\vskip 0.2 cm
Since $E(\mu^+)=-1/2$ and $E(\mu^-)=1/2$, in particular we have obtained that if $\cS_2^+\neq \emptyset$ then
$$
 \lim _{t\to +\infty}
\frac{1}{|\cS_2^+(t)|}\sum_{p\in\cS_2^+(t)} \frac{b_p}{p} =2\,
E(\mu^+)=-1\,,
$$
while for the case  $\cS_2^-\neq \emptyset$ one has
$$
\lim _{t\to
+\infty}\frac{1}{|\cS_2^-(t)|}\sum_{p\in\cS_2^-(t)} \frac{b_p}{p}
=  2\, E(\mu^-)=1\,.
$$
Hence, by Lemma \ref{expectation} the products
$$
 \prod_{p\in\cS_2^+}\frac{1}{ 1-b_p p^{- 2}+
 p^{-2}}\,, \quad  \prod_{p\in\cS_2^-}\frac{1}{ 1-b_p p^{- 2}+
 p^{-2}}
$$
do not converge. To determine the orders at $s=1$ of suitable powers of them, we consider the products
$$
 \prod_{p\in\cS_2^+}\frac{1}{ (1-b_p p^{- 2}+
 p^{-2})(1-p^{-s})}\,, \quad  \prod_{p\in\cS_2^-}\frac{1}{( 1-b_p p^{- 2}+
 p^{-2})(1+p^{-s})}
$$
and we  introduce the functions
$$
  \begin{array}{lll}
\cK^+(t) &= \displaystyle{\frac{1}{|\cS_2^+(t)|}\sum_{p\in\cS_2^+(t)}\left(\frac{ b_p}{p}+1\right)}\,, &\text{if $\cS_2^+$
is non-empty,}\\[5pt]
\cK ^-(t) &=
\displaystyle{\frac{1}{|\cS_2^-(t)|}\sum_{p\in\cS_2^-(t)}\left(\frac{
b_p}{p}-1\right)}\,, &\text{if $\cS_2^-$ is non-empty,}
 \end{array}
$$
for which we already know
$$
\lim _{t\to +\infty} \cK^+(t) =
\lim _{t\to +\infty} \cK^-(t) = 0\,.
$$
Now, we want to apply Lemma \ref{expectation} to the sequences $\{ b_p+p\}_{p\in\cS_2^+}$ and $\{ b_p-p\}_{p\in\cS_2^-}$.
If we had some control of  the rate of convergence in the above
Sato-Tate type formulas, then we would obtain conditions for  the functions $\cK^+(t)$ and $\cK^-(t)$.
 As far as we know, nothing is proved in
this sense. There exists a conjecture made
first by S. Akiyama and Y. Tanigawa in \cite{akita} for rational
elliptic curves, and then by B. Mazur in~\cite{mazur} for more
general modular forms, which is in accordance with the widely
present square root accuracy. In the particular case of elliptic curves
over $\Q$ (that is, $\dim A_f=1$), the conjecture claims that the
error term
$$
E(a,t)=\displaystyle{ \frac{|\{p\in \cP(t): \frac{a_p}{2\sqrt
p}\leq a\}|}{|\cP(t)|}-\frac{2}{\pi}\int_ {-1}^a \sqrt{1-x^2}\,
dx}\,,
$$
verifies
$$
E(a,t)=O(t^{-\frac12+\varepsilon}).
$$
The conjecture has been tested, and even refined, numerically by
W. Stein in \cite{stein} for rational elliptic curves. It is
interesting to note that his data show a relation between the rate
of convergence of $E(a,t)$ and the rank of the elliptic curve.

Nevertheless, we only need to have a control of the rate of convergence of the mathematical expectations. This is the content of the following proposition.

\begin{prop}\label{rate} Let $\cK^+(t)$, $\cK^-(t)$ as above. Then,
$$
\displaystyle{\int_{p_1}^{+\infty}\frac{\mid\cK^+(t)\mid}{t\log t}\,
dt}<+\infty \,,\qquad
\displaystyle{\int_{p_1}^{+\infty}\frac{\mid\cK^-(t)\mid}{t\log t}\,
dt}<+\infty\,.
$$
\end{prop}
{\bf Proof.}  It is well-known that the symmetric square $L$-function
$$
L((\operatorname{Symm}^2\, f)\otimes \chi,s)=\prod_{p\in\cP} \prod_{i=0}^2 (1-\alpha_p^{n-i}\beta_p ^i\,\chi(p) /p^s)^{-1},
$$
is an entire $L$-function as those considered in Section 5 of \cite{iwaniec}. In fact, Part 2 of Theorem~5.44 of~\cite{iwaniec} give us
a zero free region like
$1-\frac{c}{\log t}$, for some constant~$c$ depending on the form $f$. Now, we note that
its logarithmic derivative is given by
$$
\frac{d}{ds} \log\left(L((\operatorname{Symm}^2\, f)\otimes \chi,s)\right)=-\displaystyle{\sum_{p\in\cP}\sum_{k\geq 1} \frac{\sum_{i=0}^2(\alpha_p^{2-i}\beta_p ^i\chi(p))^k\,\log p}{p^{k\, s}}}
$$
and, by definition, we have the bound
$$\frac{\sum_{i=0}^2\alpha_p^{2-i}\beta_p ^i\,\chi(p)\,\log p}{p}\ll c\log x\,,
$$
for some constant $c$ and any $p\le x$.
Hence, the crude bound ($5.48$) in page $110$ of \cite{iwaniec}  applies, and we can use
Theorem~$5.13$ of~\cite{iwaniec} to get, in particular
$$
\sum_{p\in\cP(x)}\frac{\sum_{i=0}^2\alpha_p^{2-i}\beta_p ^i\,\chi(p)\,\log p}{p}=O\left(e^{-c\sqrt{\log x}}\right).
$$
Now we just have to apply Abel summation, to get
$$
\sum _{p\in\cP( x)}\frac{\sum_{i=0}^2 \alpha_p^{2-i}\beta_p ^i\chi(p)}{p}=O\left(e^{-K\sqrt{\log x}}\right),
$$
for some constant $K$. Summing over the characters mod $M$, and using orthogonality we get
\begin{equation}\label{eq:sym2}
\sum _{\substack{p\in\cP( x) \\ p\equiv m\,(\bmod M)}}\frac{\sum_{i=0}^2 \alpha_p^{2-i}\beta_p ^i}{p}=O\left(e^{-K\sqrt{\log x}}\right),
\end{equation}
 which clearly implies the result.
\hfill$\Box$
\begin{rem} 
It is possible to prove an analogous estimate as (\ref{eq:sym2}) for any 
symmetric $n$-th power.
From  part 2 of Theorem B in \cite{BGHT}, the function 
$L((\operatorname{Symm}^n\, f)\otimes \chi,s)$ has meromorphic continuation to the whole complex plane and is holomorphic and non-zero in $\Re(s)\geq n/2+1$.  From the proof of Theorem B, we know that  $L((\operatorname{Symm}^n\, f)\otimes \chi,s)$ is the quotient of an automorphic representation arising from RAESDC representations $\pi$ of $GL_{n+1}\left(\mathbb A_\mathbb L\right)$. Hence, we are in the right position to apply for example Theorem $5.42$ of \cite{iwaniec}  (see also \cite{sarnak}) and obtain a zero free region for $L((\operatorname{Symm}^n\, f)\otimes \chi,s)$ of the form $1-\frac{c}{\log t}$ for some constant $c$ depending on $n$. Now, we just have to consider its logarithmic derivative to obtain again the bound
$$\frac{\sum_{i=0}^n(\alpha_p^{n-i}\beta_p ^i\, \chi(p))\,\log p}{p^{n/2}}\ll c\log x\,,
$$
for some constant $c$ depending on $n$ and any $p\le x$.
Hence, the crude bound ($5.48$) in page $110$ of \cite{iwaniec}  applies, and we can use Theorem $5.13$ of the same reference to get the same estimation as in the above proof for every~$n\geq 1$.
\end{rem}

\vskip 0.2 cm

We shall need the following lemma. Recall that $\Gal (\LL/\Q)$ is
the compositum of a polyquadratic extension of $\Q$ and the cyclic
extension $\overline{\Q}^{\ker \varepsilon}$. We consider the
integers $n$, $n_1$ and $n_2$ defined as follows
$$
\displaystyle{n=\ord  \varepsilon\,,\quad
2^{n_1}=\frac{[\LL:\Q]}{n} \,,\quad  n_2= \frac{(-1)^n+1}{2}}\,.
$$
Notice that $\cS_2^+ =\emptyset$, resp. $\cS_2^-= \emptyset$, if
and only if $n_1=0$, resp. $n_2=0$.

\begin{lema}\label{Dedekind} Assume  $\cS_2\neq \emptyset$. One has,
\begin{itemize}
\item[(i)] the functions
$$
\displaystyle{G^+(s)=\left(\prod_{p\in\cS_2^+}\frac{1}{1+p^{-s}}\right)^{[\LL:\Q]}\,,\quad
G^-(s)=\left(\prod_{p\in\cS_2^+}\frac{1}{1-p^{-s}}\right)^{[\LL:\Q]}}
$$
are meromorphic on $\Re(s)>1/2$ and
$\ord_{s=1}G^+(s)=-\ord_{s=1}G^-(s)= n_1$.
\item[(ii)] The functions
$$
\displaystyle{H^+(s)=\left(\prod_{p\in\cS_2^-}\frac{1}{1+p^{-s}}\right)^{[\LL:\Q]}\,,\quad
H^-(s)=\left(\prod_{p\in\cS_2^-}\frac{1}{1-p^{-s}}\right)^{[\LL:\Q]}}
$$
are meromorphic on $\Re(s)>1/2$ and
$\ord_{s=1}H^+(s)=-\ord_{s=1}H^-(s)= n_2$.
\end{itemize}

These functions are taken to be the constant $1$ if the set of
primes in the corresponding product is empty.
\end{lema}

\noindent{\bf Proof.} Let $\M$ be a subfield of $\LL$ such that
$[\LL:\M]=2$. Let $\cS_\M$ be the set of primes of~$\cS_2$ that
split completely in $\M$. Recall that the asterisk means to
exclude Euler factors corresponding to primes dividing $N$.
 By using Dedekind zeta functions, we know that:
 \begin{itemize}
 \item The function
 $$
  S_1(s):=\left(\prod_{p\in\cS_1}\frac{1}{1-p^{-s}}\right)^{[\LL:\Q]}= \zeta_\LL^*(s)/\left(\prod_{\{p: d(p)>1\}}\frac{1}{1-p^{-s\, d(p)}}\right)^{[\LL:\Q]/d(p)}
 $$
 is meromorphic on $\Re(s)>1/2$ and its order at $s=1$ is $-1$.
 \item The function
 $$
 S_2(s):=\left(\prod_{p\in\cS_1}\frac{1}{1-p^{-s}}\prod_{p\in \cS_\M}\frac{1}{1+p^{-s}}\right)^{[\LL:\Q]/2}
 $$
  is meromorphic on $\Re(s)>1/2$ and has the same  order at $s=1$ as $ \zeta_\LL^*(s)/\zeta_{\M}^*(s)$, which
 is $0$.
 \end{itemize}
 Hence, the function
 $$
 S_2(s)^2/S_1(s)=\left(\prod_{p\in\cS_\M}\frac{1}{1+p^{-s}}\right)^{[\LL:\Q]}
 $$
 has a zero of order $1$ at $s=1$.

Now, the statement follows from the facts that  $\cS_2^+$, resp.
$\cS_2^-$, is the disjoint union of $n_1$, resp. $n_2$, sets
$\cS_{\M_i}$ such that $[\LL:\M_i]=2$ and the functions
$G^+(s)\,G^-(s)$ and $H^+(s)\,H^-(s)$ are holomorphic and non-zero
at $s=1$. \hfill$\Box$

Finally, we obtain the main theorem for the non-CM case.

\begin{teo}\label{mainnoCM}
Let $f$ be a newform without CM  such that  $\cS_2\neq \emptyset$.
Let   $n_1$ and  $n_2$ as above.  With the notations as in Section~3, then
$$
\begin{array}{lr}
\ord_{s=1} L(f/\LL,s)= &
\displaystyle{\frac{1}{2 }\left( \ord_{s=1} L(f,\cS_1,s)^{2[\LL:\Q]}+n_1-n_2\right)}\,,\\[9 pt]
\ord_{s=1} L(B_f/\LL,s)= &
\displaystyle{\frac{t}{2 }\left( \ord_{s=1} L_1(s)^{2[\LL:\Q]}+(n_1-n_2)[\F:\Q]\right)}\,,\\[9 pt]
\ord_{s=1} L(A_f/\LL,s)= &
\displaystyle{\frac{[\E:\F]}{2 } \ord_{s=1} L_1(s)^{2[\LL:\Q]}+\frac{1}{2}(n_1-n_2)[\E:\Q]}\,.
\end{array}
$$
Moreover,
if $A_f$ is isogenous to $\Res_{\LL/\Q}(B_f)$, then
$$
\ord_{s=1} L(A_f/\Q,s)=
\displaystyle{\frac{1}{2 }\left( \ord_{s=1} L_1(s)^{2[\E:\F]}+(n_1-n_2)[\F:\Q]\right)}\,.
$$

\end{teo}

\noindent{\bf Proof.} By Proposition \ref{rate} and  Lemmas~\ref{convergence},
\ref{expectation} and \ref{Dedekind}, it follows that the product
$$
\left(L(f,\cS_2,s)L(f^-,\cS_2,s)\right)^{[\LL:\Q]}  G^-(s) H^+(s)
$$
converges at $s=1$. Hence,
$$
\ord_{s=1} \left(L(f,\cS_2,s)L(f^-,\cS_2,s)\right)^{[\LL:\Q]}=\ord_{s=1} G^+(s)+\ord_{s=1} H^-(s)=n_1-n_2\,.
$$
Applying Theorem~\ref{teoremanoCM}, we obtain the
statement. \hfill$\Box$

\vskip 0.1 cm

As an application, we get the following.
\begin{prop}
Suppose that $A_f$ is an abelian surface with quaternionic
multiplication. Then,
$$
\ord_{s=1} L(A_f/\Q,s) = \frac{1}{2} \, \left( \ord_{s=1} L(f,\cS_1,s)^4 +1 \right)\,.
$$
In particular, $\ord_{s=1} L(f,\cS_1,s)^4$ is an odd integer $\geq
-1$.
\end{prop}

\noindent{\bf Proof.} We have that $\dim A_f=2$, $\F=\Q$, and
$t=2$. It is known that in this case $\varepsilon=1$ and, thus,
$n_1=1$, $n_2=0$ and $\LL$ is a quadratic field. Moreover, the
Weil restriction $\Res_{\LL/\Q}(A_f)$ is isogenous over $\Q$ to
$A_f^2$. Then, Theorem \ref{mainnoCM} gives the desired formula.~\hfill$\Box$

{\bf Acknowledgements.} We would like to thank F. Chamizo, E. Friedman, B. Gross, R. Taylor and
T. Yang for their feedback during the preparation of the manuscript. Also, we appreciate a conversation with P. Sarnak about the results in \cite{BGHT}, that helped us to finish
the proof of Proposition~\ref{rate}.
\bibliographystyle{plain}
\bibliography{modularfactor}

\noindent Josep Gonz\'alez Rovira \newline
Departament de Matem\`atica Aplicada IV \newline
Universitat Polit\`ecnica de Catalunya (UPC) \newline
Av. V\'\i ctor Balaguer, s/n. \newline
E-08800 Vilanova i la Geltr\'u \newline
email: josepg@ma4.upc.edu\newline

\noindent Jorge Jim\'enez Urroz \newline
Departament de Matem\`atica Aplicada IV \newline
Universitat Polit\`ecnica de Catalunya (UPC) \newline
Edifici C3-Campus Nord  \newline
Jordi Girona, 1-3. \newline
E-08034 Barcelona \newline
email: jjimenez@ma4.upc.edu\newline

\noindent Joan-Carles Lario Loyo \newline
Departament de Matem\`atica Aplicada II \newline
Universitat Polit\`ecnica de Catalunya (UPC) \newline
Edifici Omega - Campus Nord \newline
Jordi Girona, 1-3 \newline
E-08034 Barcelona\newline
email:  joan.carles.lario@upc.edu

\end{document}